\documentclass[10pt,reqno]{amsart}
\usepackage{amsfonts,amsmath,amssymb,amsbsy,amstext,amsthm}
\usepackage{accents,color,enumerate,enumitem,float,fullpage,verbatim}

\usepackage{url}
\usepackage{microtype}
\usepackage[colorlinks=true,hyperindex, linkcolor=magenta, pagebackref=true, citecolor=cyan]{hyperref}
\usepackage[alphabetic,lite,backrefs]{amsrefs}
\usepackage{ytableau}
\usepackage{subfiles}

\usepackage{
bm,
mathbbol,
parskip
}

\usepackage{tikz,tikz-cd}
\usetikzlibrary{positioning, matrix, shapes}
\usetikzlibrary{decorations.pathmorphing}

\usepackage{lscape}

\usepackage{titlesec}
\titleformat{\section}[block]{\large\scshape\bfseries\filcenter}{\thesection.}{1em}{}		
\titleformat{\subsection}[hang]{\large\scshape\bfseries}{\thesubsection}{1em}{}			
\titleformat{\subsubsection}[hang]{\large\scshape\bfseries}{\thesubsubsection}{1em}{}		

\usepackage{multirow}
\usepackage{array}
\usepackage{booktabs}
\newcolumntype{M}[1]{>{\centering\arraybackslash}m{#1}}
\newcolumntype{N}{@{}m{0pt}@{}}
\usepackage{diagbox}
\usepackage{cancel}

\newtheorem{lemma}{Lemma}[section]
\newtheorem{theorem}[lemma]{Theorem}
\newtheorem{prop}[lemma]{Proposition}

\theoremstyle{remark}

\newcommand{\initial}{\operatorname{in}}

\newcommand{\reg}{\operatorname{reg}}

\newcommand{\pdim}{\operatorname{pdim}}

\newcommand{\codim}{\operatorname{codim}}

\newcommand{\raj}{\operatorname{raj}}

\newcommand{\ess}{\operatorname{Ess}}
\newcommand{\rk}{\operatorname{rk}}

\newcommand{\cG}{\mathcal{G}}

\newcommand{\ZZ}{\mathbb{Z}}

\title{The MatrixSchubert Package for Macaulay2}

 \author[Ayah Almousa]{Ayah Almousa}
 \address{School of Mathematics, University of Minnesota, Minneapolis, MN, 55455}
 \email{\href{mailto:almou007@umn.edu }{almou007@umn.edu}}
 \urladdr{\url{https://sites.google.com/view/ayah-almousa}}

\author[Sean Grate]{Sean Grate}
 \address{Department of Mathematics and Statistics, Auburn University, Auburn, AL, 36849}
 \email{\href{mailto:sean.grate@auburn.edu }{sean.grate@auburn.edu}}
 \urladdr{\url{https://seangrate.com/}}

\author[Daoji Huang]{Daoji Huang}
 \address{School of Mathematics, University of Minnesota, Minneapolis, MN, 55455}
 \email{\href{mailto:huan0664@umn.edu }{huan0664@umn.edu}}
 \urladdr{\url{https://www.daojihuang.me/}}

 \author[Patricia Klein]{Patricia Klein}
 \address{Department of Mathematics, Texas A \&M University, College Station, TX, 77840}
 \email{\href{mailto:pjklein@tamu.edu }{pjklein@tamu.edu}}
 \urladdr{\url{https://patriciajklein.github.io/}}

\author[Adam LaClair]{Adam LaClair}
 \address{Department of Mathematics, Purdue University, West Lafayette, IN, 47907}
 \email{\href{mailto:alaclair@purdue.edu }{alaclair@purdue.edu}}
 \urladdr{\url{https://sites.google.com/view/adamlaclair/home}}

 \author[Yuyuan Luo]{Yuyuan Luo}
 \address{Department of Mathematics, Massachusetts Institute of Technology}
 \email{\href{mailto:lyuyuan@mit.edu }{lyuyuan@mit.edu}}
 \urladdr{\url{https://www.mit.edu/~lyuyuan/}}

 \author[Joseph McDonough]{Joseph McDonough}
 \address{School of Mathematics, University of Minnesota, Minneapolis, MN, 55455}
 \email{\href{mailto:mcdo1248@umn.edu }{mcdo1248@umn.edu}}
 \urladdr{\url{https://jmcdonough98.github.io/}}


\begin{document} 
\thanks{This work was begun at a Macaulay2 mini-school and workshop funded by NSF grant DMS-2302476 and hosted by the University of Minnesota.  We thank the NSF for its support and the University of Minnesota for its hospitality.}
 \thanks{PK is partially supported by NSF grant DMS-2246962.}
 \thanks{AL was partially supported by NSF grant DMS-2100288 and by Simons Foundation Collaboration Grant for Mathematicians \#580839.}
 
\maketitle
\vspace{-1cm}
\begin{abstract}
We introduce the \texttt{MatrixSchubert} package for the computer algebra system \textit{Macaulay2}. This package has tools to construct and study matrix Schubert varieties and alternating sign matrix (ASM) varieties. The package also introduces tools for quickly computing homological invariants of such varieties, finding the components of an ASM variety, and checking if a union of matrix Schubert varieties is an ASM variety.
\end{abstract}

\section{Introduction}

Fulton \cite{Ful92} introduced matrix Schubert varieties in the study of Schubert varieties in the complete flag variety. In that paper, Fulton showed that matrix Schubert varieties are Cohen--Macaulay and gave an attractive description of their codimension and defining equations from a combinatorial perspective.  Since that time, there has been a great deal of interest in their algebraic, geometric, and combinatorial properties (see, for example, \cites{KM05,KMY09,Hsi13,EM16,FRS16,RRR+21,PSW22,Por23,RRW23}).

More recently, Weigandt \cite{Wei17} introduced alternating sign matrix (ASM) varieites, which generalize matrix Schubert varieties, and gave a combinatorial description of their defining equations that generalizes Fulton's description in the case of matrix Schubert varieties (equivalently, when the ASM happens to be a permutation matrix).  There are to date no combinatorial descriptions of the codimension of an ASM variety nor a combinatorial criterion to determine if it is Cohen--Macaulay or even unmixed.  

The \texttt{MatrixSchubert} package implements many basic functions for permutations and ASMs (such as extending a partial ASM to an ASM, checking if one permutation avoids another, and finding the descent set of a permutation) as well as more complicated ones that rely on one or several theoretical results.  Several such examples are described in greater detail in Sections \ref{sec:schubertRegularity} and \ref{sec:schubDecomp}.

One goal of this package is to implement core results from the past three decades on commutative algebraic aspects of Schubert calculus (particularly those found in \cites{Ful92,KM05,Wei17,PSW22}) for the purpose of enjoying the theory that has already been developed.  Another is to facilitate inquiry into the many open questions that remain, including those surrounding resolutions and Betti numbers of matrix Schubert varieties and ASM varieties, Cohen--Macaulayness of ASM varieties, and codimension of ASM varieties.

\section{Background and some basic functions}

\subsection{Permutations}

Let $S_n$ denote the group of permutations of $n$ letters. For a permutation $w\in S_n$, we call the matrix that has $1$'s in the positions $(i,w(i))$ and $0$'s in all other positions the \emph{permutation matrix} of $w$.  (As a cautionary note, some authors take the convention that the matrix described here - and used throughout the package \texttt{MatrixSchubert} -  is the permutation matrix of $w^{-1}$.) 

Given $w\in S_n$, the \emph{Rothe diagram} of $w$ is \[D(w)=\{(i,j):i,j\in [n], w(i)>j, \text{ and } w^{-1}(j)>i\}.\]  The  Coxeter length of $w$ satisfies $\ell(w)= |D(w)|$. The \emph{essential set} of $w$ is 
\[\ess(w)=\{(i,j)\in D(w):(i+1,j),(i,j+1)\not \in D(w)\},\] i.e., the maximally southeast elements of the connected components of $D(w)$. 

\begin{verbatim}
i1 : needsPackage "MatrixSchubert"
i2 : w = {2,1,5,4,3};
i3 : rotheDiagram w
o3 = {(1, 1), (3, 3), (3, 4), (4, 3)}
o3 : List
i4 : essentialSet w
o4 = {(1, 1), (3, 4), (4, 3)}
o4 : List
\end{verbatim}

\subsection{Alternating Sign Matrices}\label{sec:ASMbasics}

A \emph{partial alternating sign matrix} (partial ASM) is a  matrix with entries in $\{-1,0,1\}$ so that partial sums taken along each row (and column) are all $0$ or $1$.  If the entries of each row (and column) sum to $1$ (in which case the matrix must be square), we call the partial ASM an ASM.  The ASMs whose entries all lie in $\{0,1\}$ are exactly the permutation matrices.  

 The \emph{rank function} of the $m \times n$ partial ASM $A = (A_{i,j})$ is defined by $ \rk_A(a,b)=\sum_{i=1}^a\sum_{j=1}^b A_{i,j}$ for $1 \leq a \leq m$, $1 \leq b \leq n$. (Note that if $A$ is not a permutation matrix, then $\rk_A(a,b)$ may not be the rank of the submatrix of $A$ consisting of its first $a$ rows and first $b$ columns.) 

Given a matrix $M$, let $M_{[i],[j]}$ be the submatrix of $M$ consisting of the first $i$ rows and $j$ columns.   Given an $m \times n$ partial ASM $A$, we define the \emph{ASM variety} of $A$ to be \[X_A=\{M\in \mbox{mat}(m,n): \rk(M_{[i],[j]})\leq \rk_A(i,j) \text{ for all } (i,j)\leq (m,n)\}.\]  If $A \in S_n$, we call $X_w$ a \emph{matrix Schubert variety}.  For background on matrix Schubert varieties, including a geometric motivation for their definition and a description of their connection to Schubert varieties, see \cites{Ful92,MS05}.

We call the defining radical of $X_A$ in $R = \kappa[z_{1,1},\ldots, z_{m,n}]$, $\kappa$ an arbitrary field, an \emph{ASM ideal} or, when $A \in S_n$, a \emph{Schubert determinantal ideal}.

\begin{prop}[{\cite[Proposition~3.3]{Ful92}}]
For $w \in S_n$, $I_w$ is prime, $\mbox{codim}(I_w)=\ell(w)$, $R/I_w$ is Cohen--Macaulay.
\end{prop}

The function {\tt{schubertCodim}} $w$ computes the codimension of $I_w$ (that is, the codimension of spec$(R/I_w)$ in spec$(R)$) using \cite[Proposition~3.3]{Ful92} and the equality $\ell(w) = |D(w)|$. Returning to the example $w=21543$, one may compare the outcome below to the size of $D(w)$, computed above.
\begin{verbatim}
i5 : schubertCodim w
o5 = 4
\end{verbatim}

In \cite[Lemma~3.10]{Ful92}, Fulton gave a generating set for the Schubert determinantal ideal $I_w$.  Fix an $m \times n$ generic matrix $Z=(z_{i,j})$. We write $I_k(Z_{[i],[j]})$ for the ideal of $R$ generated by the $k$-minors in the submatrix $Z_{[i],[j]}$ of $Z$ consisting of its first $i$ rows and first $j$ columns.  Then
\[I_w=\sum_{(i,j)\in\ess(w)}I_{\rk_w(i,j)+1}(Z_{[i],[j]}),\]  and we call these generators the \emph{Fulton generators}.  There is a generalization of the Fulton generators for ASM ideals (see \cite[Lemma~5.9]{Wei17}).

The function {\tt schubertDeterminantalIdeal} takes in either a permutation (as a list representing its one-line notation) or a partial ASM matrix and produces its ASM ideal via its Fulton generators (which typically do not form a minimal generating set).

\begin{verbatim}
i6 : v = {3,1,4,2};
i7 : schubertDeterminantalIdeal v
o7 = ideal (z   , z   , - z   z   + z   z  , - z   z   + z   z    , - z   z   + z   z   )
             1,1   1,2    1,2 2,1   1,1 2,2     1,2 3,1   1,1  3,2     2,2 3,1   2,1 3,2
             \end{verbatim}

A term order on $R$ so that the lead term of the determinant of any submatrix of $Z$ is the product of terms along its antidiagonal is called an \emph{antidiagonal} term order.

We now state a pivotal result, due to Knutson and Miller in the case of Schubert determinantal ideals.  One can extend the result to ASM ideals using either using Frobenius splitting \cite{Knu09} or a combinatorial argument \cite{Wei17}.  For full details on the latter, see \cite{KW}.

\begin{theorem}\cites{KM05, Knu09, Wei17, KW}\label{thm:antidiagGrob}
Fix a partial ASM $A$. The Fulton generators of $I_A$ form a Gr\"obner basis under any antidiagonal term order $<$.  Consequently, $\initial_< (I_A)$ is radical.
\end{theorem}

Rather than computing an ASM ideal and then afterwards computing a Gr\"obner basis, these theorems allow us to get an antidiagonal initial ideal directly from the permutation or ASM.

\begin{verbatim}
i8 : A = matrix{{0,0,1,0},{1,0,-1,1},{0,0,1,0},{0,1,0,0}};
i9 : antiDiagInit A
o9 = monomialIdeal (z   , z   , z   z    , z   z   , z   z   )
                     1,1   1,2   1,3  2,1  1,3  2,2  2,2  3,1
\end{verbatim}

\section{Rank tables}

The rank function $\rk_A(i,j)$ of an $m \times n$ partial ASM used in the definition ASM ideal can be applied to all $(i,j) \leq (m,n)$.  The entire function can be computed and recorded in an $m \times n$ matrix, implemented via the  method \texttt{rankTable}.  

\begin{verbatim}
i10 : M = matrix{{0,1,0},{1,-1,0}};
               2        3
o10 : Matrix ZZ  <--- ZZ
i11 : rankTable M
o11 = | 0 1 1 |
      | 1 1 1 |
               2        3
o11 : Matrix ZZ  <--- ZZ
\end{verbatim}

Given any rank table that could be constructed from a partial ASM, the method \texttt{rankTableToASM} produces the unique partial ASM of the same size as the input having that rank table.  

\begin{verbatim}
i12 : rankTableToASM matrix{{0,1,1},{0,1,1},{1,2,2}}
o12 = | 0 1 0 |
      | 0 0 0 |
      | 1 0 0 |
               3        3
o12 : Matrix ZZ  <--- ZZ
\end{verbatim}

For the typical ASM, many different rank tables could be used to construct the same ASM variety $X_A$.  The method \texttt{rankTableToASM} expects the (unique) matrix with minimum possible entries, which is the one constructed in Section \ref{sec:ASMbasics} (see \cite[Lemma 1]{RR86}).  If the user has a non-minimal rank table (as a matrix of non-negative integers), the function \texttt{rankTableFromMatrix} transforms the non-minimal rank table into a minimal rank table.

\begin{verbatim}
i13 : rankTableFromMatrix matrix{{0,1,2},{0,4,1},{8,2,4}}
o13 = | 0 1 1 |
      | 0 1 1 |
      | 1 2 2 |
               3        3
o13 : Matrix ZZ  <--- ZZ
\end{verbatim}

Among other uses, rank tables also facilitate efficient addition of ASM ideals.  Every sum of ASM ideals is again an ASM ideal \cite[Section 3]{Wei17}, and the rank table of the sum is the entrywise minimum of the rank tables of the ASMs appearing as summands.  If $A_1, \ldots, A_k$ are all $m\times n$ partial ASMs, the method \texttt{schubertAdd} computes the rank tables of the $A_i$ (via \texttt{rankTable}), takes entrywise minima (via  \texttt{entrywiseMinRankTable}), and computes an ASM ideal from that rank table. Both the ASM and its rank table are saved in the cache of the newly computed ASM ideal.

\begin{verbatim}
i14 : N = matrix{{1,0,0},{0,0,1}}
o14 = | 1 0 0 |
      | 0 0 1 |
              2        3
o14 : Matrix ZZ  <--- ZZ
i15 : idealSum = schubertAdd{M,N}
o15 = ideal (z   , - z   z    + z    z   )
              1,1     1,2  2,1   1,1  2,2
o15 : Ideal of QQ[z   ..z   ]
                   1,1   2,2
i16 : peek idealSum.cache
o16 = CacheTable{ASM => | 0 1  0 0 |     }
                        | 1 -1 0 1 |
                        | 0 1  0 0 |
                        | 0 0  1 0 |
                 rankTable => | 0 1 1 1 |
                              | 1 1 1 2 |
                              | 1 2 2 3 |
                              | 1 2 3 4 |
i17 : getASM idealSum
o17 = | 0 1  0 0 |
      | 1 -1 0 1 |
      | 0 1  0 0 |
      | 0 0  1 0 |

               4        4
o17 : Matrix ZZ  <--- ZZ
\end{verbatim}

\section{Pattern Avoidance}

The \texttt{MatrixSchubert} package has functions to test pattern avoidance for permutations. 

A permutation is called \emph{vexillary} if it avoids the permutation $2143$. The class of one-sided ladder determinantal ideals coincides exactly with the class of vexillary matrix Schubert varieties.  The vexillary condition has a large number of equivalent definitions.  We direct the reader to \cite[Section 3.2]{KMY09} for many of them. Testing if a permutation is vexillary is implemented via the \texttt{isVexillary} function. 

\begin{verbatim}
i18 : w = {7,2,5,8,1,3,6,4};
i19 : isVexillary w
o19 = false
i20 : w = {1,6,9,2,4,7,3,5,8};
i21 : isVexillary w
o21 = true
\end{verbatim}

A permutation $w$ is \emph{CDG} if it avoids all eight of the following patterns:
\begin{equation*}
    13254, 21543, 214635, 215364, 215634, 241635, 315264, 4261735.
\end{equation*}

The class of CDG permutations was named in \cite{HPW22}, where a diagonal Gr\"obner basis was conjectured for the class (proved in \cite{Kle}).  Every CDG permutation is vexillary, and the CDG permutations form the largest named class of permutations for which a diagonal Gr\"obner basis of their matrix Schubert varieties is known.  Testing if a permutation is CDG is implemented via the \texttt{isCDG} function.

\begin{verbatim}
i22 : w = {5,7,2,1,6,4,3};
i23 : isCDG w
o23 = false
i24 : w = {1,3,5,7,2,4,6};
i25 : isCDG w
o25 = true
\end{verbatim}

We say that $w$ is \emph{Cartwright--Sturmfels} if it avoids all of the following twelve patterns:
\begin{equation*}
    12543, 13254, 13524, 13542, 21543, 125364, 125634, 215364, 215634, 315264, 315624, 315642.
\end{equation*}

For background on the Cartwright--Sturmfels property of ideals in general, see \cite{CDG22}.  For a proof that the Cartwright-Sturmfels property is characterized by the given pattern avoidance condition, that Cartwright-Sturmfels Schubert determinantal ideals have an universal Gr\"obner basis, and that any initial ideal of a Cartwright-Sturmfels Schubert determinantal ideal is Cohen--Macaulay, see \cite[Theorem 4.6]{CDG22}.  Every Cartwright--Sturmfels Schubert determinantal ideal is CDG.

Testing if a permutation is Cartwright-Sturmfels is implemented via the \texttt{isCartwrightSturmfels} functions.
\begin{verbatim}
i26 : w = {3,1,2,6,5,4};
i27 : isCartwrightSturmfels w;
o27 = false
i28 : w = {6,3,5,2,1,4};
i29 : isCartwrightSturmfels w
o29 = true
\end{verbatim}

More generally, \texttt{avoidsAllPatterns} inputs a permutation and a list of patterns to avoid, and determines if the permutation avoids all of the patterns. This allows users to test conjectures related to pattern avoidance in a much more general capacity.

\section{Algorithms for Castelnuovo--Mumford regularity} \label{sec:schubertRegularity}

Castelnuovo--Mumford regularity is a fundamental invariant in commutative algebra and algebraic geometry that in a rough sense gives a measure of the complexity of a module or sheaf. 

In this package we implement (as $\mathtt{schubertRegularity}$) a purely combinatorial formula developed by Peckenik, Speyer, and Weigandt \cite[Theorem 1.2]{PSW22} for computing the Castelnuovo--Mumford regularity of $S/I_w$ for arbitrary $w \in S_n$.  We also extend the functionality of $\mathtt{schubertRegularity}$ so that it can compute the Castelnuovo--Mumford regularity of the coordinate ring associated to a partial ASM by passing to the antidiagonal initial ideal, a valid strategy in light of an important theorem of Conca and Varbaro \cite{CV20} together with Theorem \ref{thm:antidiagGrob}. See Subsection \ref{subsec:ASMReg} for a fuller explanation.

\subsection{Matrix Schubert varieties}\label{subsec:schubertReg} 
The theoretical foundation for the $\mathtt{schubertRegularity}$ function is a result of Pechenik, Speyer, and Weigandt:
\begin{theorem}[{\cite[Theorem 1.2]{PSW22}}] \label{thm_reg_via_raj}
    For $w \in S_{n}$,
    \begin{align*}
        \reg(S/I_{w}) = \raj(w) - |D(w)|
    \end{align*}
    where $\raj(w)$ is the Rajchgot index of a permutation.
\end{theorem}
We refer the reader to \cite{PSW22} for the definition of the Rajchgot index.

The function $\mathtt{schubertRegularity}$ takes either a permutation in one-line notation (that is, as a list) or a partial ASM and returns the Castelnuovo--Mumford regularity of the associated coordinate ring. The computation of the Rajchot index of a permutation invovles determining longest subsequences of the permutation subject to certain conditions. We utilize memoization in determining these longest subsequences which leads to drastic speed improvements over \textit{Macaulay2}'s built-in command $\mathtt{regularity}$. (By that, we mean that there are examples of the type we would be inclined to compute in the course of research which run reliably faster, not that we have performed any type of formal efficiency analysis.  We include a couple of illustrating examples here.) The following example run on a computer with an AMD Ryzen 5 5600U processor demonstrates the extent of this speed-up.


\begin{verbatim}
i30 : w = {1,2,3,9,8,4,5,6,7};
i31 : I = antiDiagInit(w, CoefficientRing=>ZZ/3001);
                        ZZ
o31 : MonomialIdeal of ----[z   ..z   ]
                       3001  1,1   9,9
i32 : M = comodule I;
i33 : time regularity M
     -- used 1.41504 seconds
o33 = 6
i34 : time schubertRegularity w
     -- used 0.000181139 seconds
o34 = 6
i35 : time schubertRegularity random toList (1 .. 100)
     -- used 4.03299 seconds
o35 = 1925
\end{verbatim}

\subsection{ASM Varieties}\label{subsec:ASMReg} The function $\mathtt{schubertRegularity}$ also accepts as input a partial ASM. First, $\mathtt{schubertRegularity}$ checks whether the matrix is a permutation matrix; in which case the Castelnuovo--Mumford regularity is computed via Theorem \ref{thm_reg_via_raj}. Otherwise, the antidiagonal initial ideal of the ASM is computed, and the built-in {\tt Macaulay2} command $\mathtt{regularity}$ is used. The Castelnuovo--Mumford regularity of the quotient by the initial ideal will coincide with the Castelnuovo--Mumford regularity of the coordinate ring corresponding to the ASM by {\cite[Corollary 2.7]{CV20}} since the antidiagonal initial ideal is squarefree \cites{KM05,Wei17}.
\begin{verbatim}
i36 : A = matrix{{0,0,1,0},{0,1,-1,1},{1,-1,1,0},{0,1,0,0}};
              4        4
o36 : Matrix ZZ  <--- ZZ
i37 : time regularity comodule schubertDeterminantalIdeal A
     -- used 0.00968312 seconds
o37 = 1
i38 : time schubertRegularity A
     -- used 0.0100184 seconds

o38 = 1
i39 : B = matrix{{1,0,0,0,0,0,0,0},{0,1,0,0,0,0,0,0},{0,0,0,0,1,0,0,0},
         {0,0,0,0,0,1,0,0},{0,0,1,0,0,-1,1,0},{0,0,0,1,-1,1,0,0},
         {0,0,0,0,1,0,0,0},{0,0,0,0,0,0,0,1}};
              8        8
o39 : Matrix ZZ  <--- ZZ
i40 : time regularity comodule schubertDeterminantalIdeal B
     -- used 1.01169 seconds
o40 = 8
i41 : time schubertRegularity B
     -- used 0.08511 seconds
o41 = 8
\end{verbatim}

We are a little bit disappointed that, as we see above with the ASM $A$, the command \texttt{schubertRegularity} is sometimes slower than the already-available \texttt{regularity} command applied to an ASM ideal.  However, as we see above with the ASM $B$, \texttt{schubertRegularity} can be a meaningful improvement over \texttt{regularity} in other cases.\footnote{We expect that the issue is that we are for some reason describing the sets of rows and columns whose minors define the ASM ideal (or its antidiagonal initial ideal) less efficiently than the \texttt{minors} command does.  This explanation is compatible with our experience that \texttt{antiDiagInit} is typically slower than \texttt{schubertDeterminantalIdeal} for dominant permutations, i.e., permutations indexing Schubert determinantal ideals that are already monomial ideals.} As a rule of thumb, we recommend \texttt{schubertRegularity} for ASM ideals with more generators of higher degrees and \texttt{regularity} for ASM ideals with fewer generators or generators in lower degrees.

\section{Schubert and Grothendieck polynomials}
The \texttt{MatrixSchubert} package provides functions to compute Schubert, double Schubert, and Grothendieck polynomials for permutations. We give a brief overview of how these families of polynomials are constructed, and refer the reader to \cite{KM05} for a more detailed (and somewhat more general) treatment. Let $\mathbf{x} = \{x_1, \ldots, x_n\}$, and let $S$ denote the polynomial ring in $\mathbf{x}$ over the field $\kappa$.  Let $\partial_i$ be the $i$'th divided difference operator, which sends $f \in S$ to
\[
\partial_i(f) = \frac{f(x_1,\ldots,x_n) - f(x_1,..,x_{i+1},x_i,..,x_n)}{x_i-x_{i+1}}. 
\]
It is not obvious that $\partial_i(f) \in S$, but it is true.  Let $w_0 = nn-1\cdots 321$, i.e., the longest word in $S_n$.  For $w \in S_n$, the \emph{Schubert polynomial} $\mathfrak{S}_w$ is defined recursively as follows:
\[
\mathfrak{S}_{w_0}(\mathbf{x}) = \prod_{i=1}^{n}x_i^{n-i} \in \ZZ[\mathbf{x}] \qquad \text{and}\qquad \mathfrak{S}_{ws_i}(\mathbf{x}) = \partial_i \mathfrak{S}_w(\mathbf{x})
\]
where $s_i$ is a right descent of $w$ (i.e., $w(i)>w(i+1)$). The \emph{double Schubert polynomial} $\mathfrak{S}_w(\mathbf{x},\mathbf{y})$ is defined using the same recursion, but with initial condition $\mathfrak{S}_{w_0}(\mathbf{x},\mathbf{y}) = \prod_{i+j\leq n}(x_i-y_j) \in \ZZ[\mathbf{y}][\mathbf{x}]$.
Finally, the \emph{Grothendieck polynomial} $\cG_w(\mathbf{x})$ is defined by the recursion 
\[
\cG_{w_0}(\mathbf{x}) = \prod_{i=1}^n x_i^{n-i} \qquad\text{and}\qquad \cG_{ws_i}(\mathbf{x}) = \partial_i(\cG_w(\mathbf{x})-x_{i+1}\cG_{w}(\mathbf{x}))
\]
where $s_i$ is a right descent of $w$.
Computing these polynomials for a permutation $w \in S_n$ given in one-line notation is implemented via the \texttt{schubertPolynomial}, \texttt{doubleSchubertPolynomialnomial}, and \texttt{grothendieckPolynomial} functions respectively.
\begin{verbatim}
i42 : w = {2,1,4,3};
i43 : schubertPolynomial w
       2
o43 = x  + x x  + x x
       1    1 2    1 3
o43 : QQ[x ..x ]
          1   4
i44 : doubleSchubertPolynomialnomial w
       2                                        2
o44 = x  + x x  + x x  - 2x y  - x y  - x y  + y  - x y  + y y  - x y  + y y
       1    1 2    1 3     1 1    2 1    3 1    1    1 2    1 2    1 3    1 3
o44 : QQ[x ..x , y ..y ]
          1   4   1   4
i45 : grothendieckPolynomial w
\end{verbatim}
\smallskip
\begin{verbatim}
       2        2      2               2
o45 = x x x  - x x  - x x  - x x x  + x  + x x  + x x
      1 2 3    1 2    1 3    1 2 3    1    1 2    1 3
o45 : QQ[x ..x ]
          1   4
\end{verbatim}

The default options for computing Schubert, double Schubert, and Grothendieck polynomials use the definitions by divided difference operators, where we deterministically pick one reduced word for each $w$ to apply the divided difference operators. For Schubert polynomials, we also provide the option {\tt  Algorithm=>"Transition"}
that computes Schubert polynomials via transition equations; see, e.g., \cite{W21}. For Grothendieck polynomials, we provide three different implementations: {\tt "DividedDifference"}, {\tt "Degree"}, and {\tt "PipeDream"}. The {\tt "Degree"} option computes the twisted $K$-polynomials of the matrix Schubert variety and should not be used for any practical implementation. The {\tt "PipeDream"} option computes Grothendieck polynomials by the pipe dream formula. We provide these different options for any users who are interested in comparing efficiency of the different algorithms.

\section{Studying ASM varieties via initial ideals}\label{sec:schubDecomp}

Rank functions induce a lattice structure on the set of $n \times n$ ASMs defined by $A\geq B$  if and only if  $\rk_A(i,j)\leq \rk_B(i,j)$ for all $i,j\in [n]$.  The restriction of this partial order to $S_n$ recovers (strong) Bruhat order on $S_n$. Define \[
\text{perm}(A)=\{w\in S_n:w\geq A, \text{ and, if } w\geq v\geq A \mbox{ for some $v \in S_n$, then }  w=v\}.
\] 

\begin{prop}\cite[Proposition~5.4]{Wei17},\cite[Lemma~2.6]{KW} If $A$ is an ASM and $<$ is an antidiagonal term order, then \[
I_A = \bigcap_{w \in \text{perm}(A)} I_w \mbox{ \hspace{2cm} and \hspace{2cm}} \initial_< (I_A) = \bigcap_{w \in \text{perm}(A)} \initial_< (I_w). 
\]
\end{prop}

By combining Knutson and Miller's \cite[Theorem B]{KM05} with Bergeron and Billey's \cite[Theorem 3.7]{BB93}, one may construct a reduced word for $w \in S_n$ from the indices of the variables generating any minimal prime of $\initial_<(I_w)$.  

It is of independent combinatorial interest to understand the lattice of ASMs.  The function {\tt permSetOfASM} takes in an ASM $A$ and computes $\text{perm}(A)$ by decomposing the antidiagonal initial ideal of $I_A$ and reading a reduced word from each of the primes appearing in the decomposition.  The set of distinct permutations encountered comprises $\text{perm}(A)$ and therefore also indexes the components in a prime decomposition of $I_A$.  

\begin{verbatim}
i46 : A = matrix{{0,1,0},{1,-1,1},{0,1,0}};
i47 : permSetOfASM A
o47 = {{3, 1, 2}, {2, 3, 1}}
\end{verbatim}

The function {\tt schubertDecompose} takes in an ideal, computes its  initial ideal by the default term order in {\tt Macaulay2} (which is antidiagonal) and, from the minimal primes of that ideal, finds and returns the set of permutations with at least one reduced word given by the set of generators of one of those minimal primes.  The primary use of this function is on an ideal the user knows to be an ASM ideal $I_A$ (such as one arising directly from the matrix $A$ or as a sum of other ASM ideals), in which case the output will be $\text{perm}(A)$.

\begin{verbatim}
i48 : schubertDecompose schubertDeterminantalIdeal A
o48 = {{3, 1, 2}, {2, 3, 1}}
\end{verbatim}

If the user has an ideal and is unsure if that ideal is an ASM ideal, they may use the function {\tt isASMIdeal}, which takes in an ideal $I$ and first applies {\tt schubertDecompose}.  It then takes entrywise maxima (using {\tt entrywiseMaxRankTable}) among the rank tables of the permutations found from {\tt schubertDecompose}.  Using {\tt rankTableToASM}, it constructs the partial ASM $A$ whose ASM ideal is determined by that rank table.  Finally, {\tt isASMIdeal} returns a boolean indicating if $I_A$ is equal to the input ideal.  In case it is, {\tt isASMIdeal} caches the partial ASM $A$ so that $I=I_A$.  The ASM $A$ may then be retrieved via {\tt getASM}.

\begin{verbatim}
i49 : I1 = schubertDeterminantalIdeal {3,4,1,2};
o49 : Ideal of QQ[z   ..z   ]
                   1,1   4,4
i50 : I2 = sub(schubertDeterminantalIdeal {3,2,4,1},ring I1);
o50 : Ideal of QQ[z   ..z   ]
                   1,1   4,4
i51 : I = intersect(I1,I2);
o51 : Ideal of QQ[z   ..z   ]
                   1,1   4,4
i52 : isASMIdeal I
o52 = true
i53 : A = getASM I
o53 = | 0 0  1 0 |
      | 0 1  0 0 |
      | 1 -1 0 1 |
      | 0 1  0 0 |
\end{verbatim}

The antidiagonal initial ideal of an ASM ideal can be used for additional computations in light of recent and impactful results of Conca and Varbaro \cite{CV20}. Fix a partial ASM $A$ and antidiagonal term order $<$.  A result from \cite{CV20} closely related to that discussed in the section on Castelnuovo--Mumford regulartiy states that $R/I_A$ is Cohen--Macaulay if and only if $R/\initial_< (I_A)$ is Cohen--Macaulay (which again uses that $\initial_< (I_A)$ is radical).
The function {\tt isSchubertCM} assesses Cohen--Macaulayness of $R/I_A$ by checking $\pdim(\initial_<(I_A))==\codim(\initial_<(I_A))$.

We continue from the example above and then consider a familiar non-Cohen--Macaulay variety.

\begin{verbatim}
i54 : isSchubertCM A
o54 = true
i55 : B = matrix{{0,0,1,0,0},{0,0,0,1,0},{1,0,-1,0,1},{0,1,0,0,0},{0,0,1,0,0}};
                 5        5
o55 : Matrix ZZ  <--- ZZ
i56 : trim schubertDeterminantalIdeal B
o56 = ideal (z   , z   , z   , z   , z   z   , z   z   , z   z   , z   z   )
              2,2   2,1   1,2   1,1   2,3 3,2   1,3 3,2   2,3 3,1   1,3 3,1
i57 : isSchubertCM B
o57 = false
\end{verbatim}

The Stanley--Reisner complexes of antidiagonal initial ideals of ASM ideals are of independent interest.  Knutson and Miller \cite{KM05} introduced subword complexes, of which Stanley--Reisner complexes of antidiagonal initial ideals of Schubert determinantal ideals are the motivating example.  Given a permutation $w$, the method \texttt{subwordComplex} produces the Stanley--Reisner complex of the antidiagonal initial ideal of $I_w$.

\begin{verbatim}
i58 : w = {2,1,4,3};

i59 : netList facets subwordComplex w

       +--------------------------------------------------------+
o59 =  |z   z   z   z   z   z   z   z   z   z   z   z   z   z   |
       | 1,2 1,4 2,1 2,2 2,3 2,4 3,1 3,2 3,3 3,4 4,1 4,2 4,3 4,4|
       +--------------------------------------------------------+
       |z   z   z   z   z   z   z   z   z   z   z   z   z   z   |
       | 1,2 1,3 1,4 2,1 2,3 2,4 3,1 3,2 3,3 3,4 4,1 4,2 4,3 4,4|
       +--------------------------------------------------------+
       |z   z   z   z   z   z   z   z   z   z   z   z   z   z   |
       | 1,2 1,3 1,4 2,1 2,2 2,3 2,4 3,2 3,3 3,4 4,1 4,2 4,3 4,4|
       +--------------------------------------------------------+
i60 : v = {2,1,6,3,5,4};
i61 : # facets subwordComplex v
o61 = 35
\end{verbatim}

Moreover, Knutson and Miller \cite{KM05} showed that the prime components of the antidiagonal initial ideal of Schubert determinantal ideals are indexed by combinatorial objects called \emph{pipe dreams} (by \cite{KM05}, now adopted as standard terminology in the literature) or \emph{RC-graphs} (by \cite{BB93}, to which we refer the reader for background).  

In particular, to read off an associated prime of the antidiagonal initial ideal from a pipe dream, one simply needs to read off the locations of the $+$ tiles in the pipe dream.
This package provides the class \texttt{PipeDream} to display and manipulate pipe dreams.

\begin{verbatim}

i62 : u = {2,1,4,3,6,5};   

i63 : (pipeDreams u)_0

o63 = +/+/+/
      //////
      //////
      //////
      //////
      //////

o63 : PipeDream

i64 : (decompose antiDiagInit u)_0

o64 = monomialIdeal (z   , z   , z   )
                     1,1   1,3   1,5

o64 : MonomialIdeal of QQ[z   ..z   ]
                          1,1   6,6
\end{verbatim}

The generator $z_{1,1}$ corresponds to the $+$ tile in location $(1,1)$ of the given pipe dream, the generator $z_{1,3}$ to the $+$ tile in location $(1,3)$, and the generator $z_{1,5}$ to the $+$ tile in location $(1,5)$.

To compare the Macaulay2 drawing of a pipe dream with those in \cite{KM05} and subsequent literature,  \texttt{+} is interpreted as a cross tile, and \texttt{/} is interpreted as an elbow tile.

\begin{figure}[h]
    \centering
    \includegraphics[width=1in]{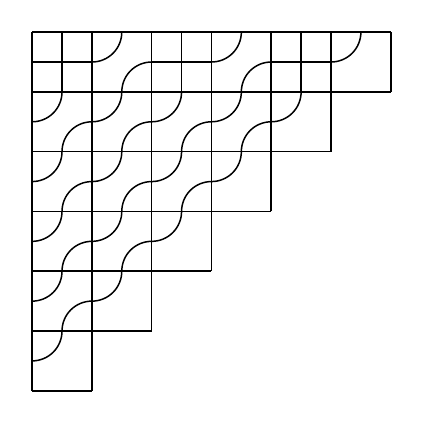}
    \caption{One pipe dream for $w = 214365$. \texttt{+} is interpreted as a cross tile, and \texttt{/} is interpreted as an elbow tile, and entries below the main antidiagonal are ignored, as is standard.}
    \label{fig:pipedreamexample}
\end{figure}

\bigskip

\section{Some diagonal term orders}
Knutson and Miller's \cite{KM05} result that every antidiagonal term order determines the same initial ideal of a given  Schubert determinantal ideal has been used heavily throughout the \verb|MatrixSchubert| package.  By contrast, different diagonal term orders can yield different initial ideals.  The various diagonal initial ideals of Schubert determinantal ideals are a topic of active research.  Given a permutation in one-line notation or a partial ASM, the functions \verb|diagLexInitSE|, \verb|diagLexInitNW|, \verb|diagRevLexInit| will each produce an initial ideal under a distinct diagonal term order.  

Precisely, \verb|diagLexInitSE| is the lexicographic order for which $z_{n,n}$ is largest and the remaining variables are ordered by reading left across the bottom row, then right to left across row $n-1$, and so on until arriving finally at $z_{1,1}$, the smallest variable.  \verb|diagLexInitNW| is the lexicographic order for which $z_{1,1}$ is largest and the remaining variables are ordered by reading right across the top row, then left to right across row $2$, and so on until arriving finally at $z_{n,n}$, the smallest variable.  And \verb|diagRevLexInit| is the reverse lexicographic order where $z_{n,1}$ is the smallest (or most penalized) variable followed by the variables encountered reading left to right along the bottom row, then left to right along row $n-1$ and so on until arriving at $z_{1,n}$.

The example $w = 214365$, taken from \cite{KW}, is the smallest example of which the authors are aware of a permutation with different initial ideals for different diagonal term orders.  The authors are unaware of any examples for which \verb|diagLexInitSE| and \verb|diagRevLexInit| produce different initial ideals.  In the example below, the Macaulay2 output recording the ambient rings of the initial ideals has been omitted for brevity. 

\begin{verbatim}
i65 : w = {2,1,4,3,6,5}
o65 = {2, 1, 4, 3, 6, 5}
o65 : List
i66 : diagLexInitSE w   
                                                                 2
o66 = monomialIdeal (z   z   z   z   z   z   , z   z   z   z   z   , 
                      5,5 4,3 3,4 3,2 2,1 1,3   5,5 4,3 3,4 2,1 1,2   
                                        z   z   z   , z   z   z   z   z   z   , z   )
                                        3,3 2,1 1,2   5,5 4,3 3,4 3,1 2,3 1,2   1,1
i67 : diagLexInitNW w
                                                                                             
o67 = monomialIdeal (z   , z   z   z   , z   z   z   z   z   ,  
                      1,1   1,2 2,1 3,3   1,2 2,1 3,4 4,3 5,5   
                                  2
                                z   z   z   z   z   z   , z   z   z   z   z   z   )
                                 1,2 2,3 3,1 3,4 4,3 5,5   1,3 2,1 3,2 3,4 4,3 5,5
i68 : diagRevLexInit w                    
                                                                 2
o68 = monomialIdeal (z   z   z   z   z   z   , z   z   z   z   z   , 
                      5,5 4,3 3,4 3,2 2,1 1,3   5,5 4,3 3,4 2,1 1,2                         
                                        z   z   z   , z   z   z   z   z   z   , z   )
                                        3,3 2,1 1,2   5,5 4,3 3,4 3,1 2,3 1,2   1,1
\end{verbatim}

\section{Available Examples of ASMs}
For the convenience of the user, we provide a complete list of ASMs up to ASM(7). It can be accessed as 
\verb|ASMFullList n| for any $1 \leq n \leq 7$, which returns a list of objects of type \verb|Matrix|. To access a list of $m$ random examples of ASMs of size $n$, \verb|ASMRandomList(n,m)|
returns a random length-$m$ list of ASMs which are size $n$, presented as \verb|Matrix| objects. 

\begin{verbatim}
i69 : ASMRandomList(5,4)

o69 = {| 1 0 0 0  0 |, | 0 0 0 1  0 |, | 0 0  0 1  0 |, | 0 0 0  0 1 |}
       | 0 0 1 0  0 |  | 1 0 0 -1 1 |  | 0 1  0 -1 1 |  | 0 0 1  0 0 |
       | 0 0 0 1  0 |  | 0 1 0 0  0 |  | 1 -1 1 0  0 |  | 1 0 0  0 0 |
       | 0 1 0 -1 1 |  | 0 0 0 1  0 |  | 0 0  0 1  0 |  | 0 1 -1 1 0 |
       | 0 0 0 1  0 |  | 0 0 1 0  0 |  | 0 1  0 0  0 |  | 0 0 1  0 0 |
\end{verbatim}

Additionally, for $n \leq 6$, lists of non-permutation ASMs that define an arithmetically Cohen--Macaulay variety, ASMs that do not define an arithmetically Cohen--Macaulay variety, and antidiagonal initial ideals of ASMs are provided, and can be accessed with \verb|cohenMacaulayASMsList n|, \verb|nonCohenMacaulayASMsList n|, and \verb|initialIdealsList n|, respectively.

It is well known that there are $429$ $5 \times 5$ ASMs, of which $5!$ are permutation matrices.  Each non-permutation $5 \times 5$ ASM is on exactly one of the lists \verb|cohenMacaulayASMsList 5| or \verb|nonCohenMacaulayASMsList 5|.

\begin{verbatim}
i70 : CM = cohenMacaulayASMsList 5;
i71 : NCM =  nonCohenMacaulayASMsList 5;
i72 : #CM+#NCM+5! == 429
o72 = true
i73 : #ASMFullList 5
o73 = 429
\end{verbatim}

\section{Acknowledgements}
The authors thank Shiliang Gao, Pooja Joshi, and Antsa Tantely Fandresena Rakotondrafara, contributors to the MatrixSchubert package who provided code, improved existing code, or enhanced documentation.  They also thank Anton Leykin, Mike Stillman and Gregory Smith, who answered countless questions both during the Macaulay2 workshop, where this package was begun, and after. They thank Mahrud Sayrafi for his helpful code review and for implementing the PipeDream class. The authors thank Anna Weigandt and Jonathan Monta{\~n}o for helpful conversations.
Finally, the authors would also like to thank Christine Berkesch, Michael Perlman, and Mahrud Sayrafi for co-organizing the workshop with the first author where this project began.

\bibliographystyle{amsplain}
\bibliography{biblio}

\begin{bibdiv}
\begin{biblist}

\bib{BB93}{article}{
      author={Bergeron, Nantel},
      author={Billey, Sara},
       title={R{C}-graphs and {S}chubert polynomials},
        date={1993},
        ISSN={1058-6458},
     journal={Experiment. Math.},
      volume={2},
      number={4},
       pages={257\ndash 269},
  url={http://projecteuclid.org.srv-proxy1.library.tamu.edu/euclid.em/1048516036},
      review={\MR{1281474}},
}

\bib{CDG22}{article}{
      author={Conca, Aldo},
      author={De~Negri, Emanuela},
      author={Gorla, Elisa},
       title={Radical generic initial ideals},
        date={2022},
        ISSN={2305-221X},
     journal={Vietnam J. Math.},
      volume={50},
      number={3},
       pages={807\ndash 827},
  url={https://doi-org.srv-proxy1.library.tamu.edu/10.1007/s10013-022-00551-w},
      review={\MR{4447407}},
}

\bib{CV20}{article}{
      author={Conca, Aldo},
      author={Varbaro, Matteo},
       title={Square-free {G}r\"{o}bner degenerations},
        date={2020},
        ISSN={0020-9910},
     journal={Invent. Math.},
      volume={221},
      number={3},
       pages={713\ndash 730},
  url={https://doi-org.srv-proxy1.library.tamu.edu/10.1007/s00222-020-00958-7},
      review={\MR{4132955}},
}

\bib{EM16}{article}{
      author={Escobar, Laura},
      author={M\'{e}sz\'{a}ros, Karola},
       title={Toric matrix {S}chubert varieties and their polytopes},
        date={2016},
        ISSN={0002-9939,1088-6826},
     journal={Proc. Amer. Math. Soc.},
      volume={144},
      number={12},
       pages={5081\ndash 5096},
         url={https://doi.org/10.1090/proc/13152},
      review={\MR{3556254}},
}

\bib{FRS16}{article}{
      author={Fink, Alex},
      author={Rajchgot, Jenna},
      author={Sullivant, Seth},
       title={Matrix {S}chubert varieties and {G}aussian conditional
  independence models},
        date={2016},
        ISSN={0925-9899,1572-9192},
     journal={J. Algebraic Combin.},
      volume={44},
      number={4},
       pages={1009\ndash 1046},
         url={https://doi.org/10.1007/s10801-016-0698-2},
      review={\MR{3566228}},
}

\bib{Ful92}{article}{
      author={Fulton, William},
       title={Flags, {S}chubert polynomials, degeneracy loci, and determinantal
  formulas},
        date={1992},
        ISSN={0012-7094},
     journal={Duke Math. J.},
      volume={65},
      number={3},
       pages={381\ndash 420},
  url={https://doi-org.srv-proxy1.library.tamu.edu/10.1215/S0012-7094-92-06516-1},
      review={\MR{1154177}},
}

\bib{HPW22}{article}{
      author={Hamaker, Zachary},
      author={Pechenik, Oliver},
      author={Weigandt, Anna},
       title={Gr\"{o}bner geometry of {S}chubert polynomials through ice},
        date={2022},
        ISSN={0001-8708,1090-2082},
     journal={Adv. Math.},
      volume={398},
       pages={Paper No. 108228, 29},
         url={https://doi.org/10.1016/j.aim.2022.108228},
      review={\MR{4379747}},
}

\bib{Hsi13}{article}{
      author={Hsiao, Jen-Chieh},
       title={On the {$F$}-rationality and cohomological properties of matrix
  {S}chubert varieties},
        date={2013},
        ISSN={0019-2082,1945-6581},
     journal={Illinois J. Math.},
      volume={57},
      number={1},
       pages={1\ndash 15},
         url={http://projecteuclid.org/euclid.ijm/1403534482},
      review={\MR{3224557}},
}

\bib{Kle}{article}{
      author={Klein, Patricia},
       title={Diagonal degenerations of matrix {S}chubert varieties},
        date={2023},
        ISSN={2589-5486},
     journal={Algebr. Comb.},
      volume={6},
      number={4},
       pages={1073\ndash 1094},
      review={\MR{4635090}},
}

\bib{KM05}{article}{
      author={Knutson, Allen},
      author={Miller, Ezra},
       title={Gr\"{o}bner geometry of {S}chubert polynomials},
        date={2005},
        ISSN={0003-486X},
     journal={Ann. of Math. (2)},
      volume={161},
      number={3},
       pages={1245\ndash 1318},
  url={https://doi-org.srv-proxy1.library.tamu.edu/10.4007/annals.2005.161.1245},
      review={\MR{2180402}},
}

\bib{KMY09}{article}{
      author={Knutson, Allen},
      author={Miller, Ezra},
      author={Yong, Alexander},
       title={Gr\"{o}bner geometry of vertex decompositions and of flagged
  tableaux},
        date={2009},
        ISSN={0075-4102,1435-5345},
     journal={J. Reine Angew. Math.},
      volume={630},
       pages={1\ndash 31},
         url={https://doi.org/10.1515/CRELLE.2009.033},
      review={\MR{2526784}},
}

\bib{Knu09}{article}{
      author={Knutson, Allen},
       title={Frobenius splitting, point-counting, and degeneration},
        date={2009},
     journal={Preprint},
        note={https://arxiv.org/abs/0911.4941},
}

\bib{KW}{article}{
      author={Klein, Patricia},
      author={Weigandt, Anna},
       title={Bumpless pipe dreams encode {G}r\"{o}bner geometry of {S}chubert
  polynomials},
        date={2022},
        ISSN={1286-4889},
     journal={S\'{e}m. Lothar. Combin.},
      volume={86B},
       pages={Art. 84, 12},
      review={\MR{4490925}},
}

\bib{MS05}{book}{
      author={Miller, Ezra},
      author={Sturmfels, Bernd},
       title={Combinatorial commutative algebra},
      series={Graduate Texts in Mathematics},
   publisher={Springer-Verlag, New York},
        date={2005},
      volume={227},
        ISBN={0-387-22356-8},
      review={\MR{2110098}},
}

\bib{Por23}{article}{
      author={Portakal, Irem},
       title={Rigid toric matrix {S}chubert varieties},
        date={2023},
        ISSN={0925-9899,1572-9192},
     journal={J. Algebraic Combin.},
      volume={57},
      number={4},
       pages={1265\ndash 1283},
         url={https://doi.org/10.1007/s10801-023-01229-3},
      review={\MR{4588133}},
}

\bib{PSW22}{article}{
      author={Pechenik, Oliver},
      author={Speyer, David~E.},
      author={Weigandt, Anna},
       title={Regularity of matrix {S}chubert varieties},
        date={2022},
     journal={S\'{e}m. Lothar. Combin.},
      volume={86B},
       pages={Art. 47, 12},
      review={\MR{4490888}},
}

\bib{RR86}{article}{
      author={Robbins, David~P.},
      author={Rumsey, Howard, Jr.},
       title={Determinants and alternating sign matrices},
        date={1986},
        ISSN={0001-8708},
     journal={Adv. in Math.},
      volume={62},
      number={2},
       pages={169\ndash 184},
         url={https://doi.org/10.1016/0001-8708(86)90099-X},
      review={\MR{865837}},
}

\bib{RRR+21}{article}{
      author={Rajchgot, Jenna},
      author={Ren, Yi},
      author={Robichaux, Colleen},
      author={St.~Dizier, Avery},
      author={Weigandt, Anna},
       title={Degrees of symmetric {G}rothendieck polynomials and
  {C}astelnuovo-{M}umford regularity},
        date={2021},
        ISSN={0002-9939,1088-6826},
     journal={Proc. Amer. Math. Soc.},
      volume={149},
      number={4},
       pages={1405\ndash 1416},
         url={https://doi.org/10.1090/proc/15294},
      review={\MR{4242300}},
}

\bib{RRW23}{article}{
      author={Rajchgot, Jenna},
      author={Robichaux, Colleen},
      author={Weigandt, Anna},
       title={Castelnuovo-{M}umford regularity of ladder determinantal
  varieties and patches of {G}rassmannian {S}chubert varieties},
        date={2023},
        ISSN={0021-8693,1090-266X},
     journal={J. Algebra},
      volume={617},
       pages={160\ndash 191},
         url={https://doi.org/10.1016/j.jalgebra.2022.11.001},
      review={\MR{4513784}},
}

\bib{Wei17}{article}{
      author={Weigandt, Anna},
       title={Prism tableaux for alternating sign matrix varieties},
        date={2017},
     journal={Preprint},
        note={arXiv:1708.07236},
}

\bib{W21}{article}{
      author={Weigandt, Anna},
       title={Bumpless pipe dreams and alternating sign matrices},
        date={2021},
     journal={Journal of Combinatorial Theory, Series A},
      volume={182},
       pages={105470},
}

\end{biblist}
\end{bibdiv}
\addcontentsline{toc}{section}{Bibliography}

\end{document}